\documentstyle[amssymb,color,epic,eepic,12pt]{article}

\newlength{\minitwocolumn}
\newcommand{\Z}{{\Bbb Z}} %??
\newcommand{\C}{{\Bbb C}} %??
\newcommand{\cL}{{\cal L}}

\newcommand{\la}{\lambda}

\newcommand{\nn}{{\nonumber}}
\newcommand{\eqref}[1]{(\ref{#1})}
\newcommand{\bea}{\begin{eqnarray}}
\newcommand{\ena}{\end{eqnarray}}
\newcommand{\beit}{\begin{itemize}}
\newcommand{\enit}{\end{itemize}}

\newcommand{\be}{\begin{eqnarray*}}
\newcommand{\en}{\end{eqnarray*}}
\newcommand{\lb}[1]{\label{#1}}
\newcommand{\End}{{\rm End}}

\newcommand{\id}{\hbox{id}}

\newcommand{\BW}[5]
{\left(\left.\matrix{#1 & #2 \cr #3 & #4 \cr}
\right| #5\right)}

\def\infq4p#1{{(#1;q^4,p)_\infty}}

\newcommand{\al}{\alpha}

\newcommand{\vep}{\varepsilon}

\font\teneufm=eufm10
\font\seveneufm=eufm7
\font\fiveeufm=eufm5
\newfam\eufmfam
\textfont\eufmfam=\teneufm
\scriptfont\eufmfam=\seveneufm
\scriptscriptfont\eufmfam=\fiveeufm

\let\goth\frak
\newcommand{\slth}{\widehat{\goth{sl}}_2}

\newcommand{\Aqp}{{\cal A}_{q,p}}

\newcommand{\Bqla}{{{\cal B}_{q,\lambda}}}

\makeatletter
\@addtoreset{equation}{section}
\makeatother

\newtheorem{thm}{Theorem}[section]
\newtheorem{prop}[thm]{Proposition}

\newtheorem{cor}[thm]{Corollary}

\newcommand{\qed}{\hspace{\fill}$\square$}
\newcommand{\vt}[3]{\vartheta_{#1}\left(\left. \frac{#2}{r} \right| {#3} \right)}
\newcommand{\vtf}[3]{\vartheta_{#1}\left(\left. {#2} \right| {#3} \right)}

\newcommand{\tp}{{\tilde{p}}}

\newcommand{\baj}{{\bar{j}}}
\newcommand{\bao}{{\bar{1}}}

\newcommand{\bal}{{\bar{l}}}

\newcommand{\nc}{\newcommand}
\nc{\sli}{\sum\limits} 
\nc{\mref}[1]{(\ref{#1})}
\nc{\vo}{v_{\gL_0}}
\nc{\vot}{v_{\gL_1+\gL_0}}
\nc{\vw}{v_{\gL_1}}
\nc{\ppmm}{\genfrac{}{}{-10pt}{10pt}{++}{--}}
\nc{\wom}[5]{\Omega\left(\left.\begin{array}{ll}{#1}&{#2}\\{#3}&{#4}\end{array}\right|{#5}\right)}
\nc{\com}[5]{\chi\left(\left.\begin{array}{ll}{#1}&{#2}\\{#3}&{#4}\end{array}\right|{#5}\right)}
\nc{\we}[5]{W\left(\left.\begin{array}{ll}{#1}&{#2}\\{#3}&{#4}\end{array}\right|{#5}\right)}
\nc{\lmat}[6]{\ell_{#6}\left(\left.\begin{array}{ll}{#1}&{#2}\\{#3}&{#4}\end{array}\right|{#5}\right)}
\nc{\lmats}[5]{L\left(\left.\begin{array}{ll}{#1}&{#2}\\{#3}&{#4}\end{array}\right|{#5}\right)}
\nc{\hmat}[6]{h_{#6}\left(\left.\begin{array}{ll}{#1}&{#2}\\{#3}&{#4}\end{array}\right|{#5}\right)}
\nc{\hmats}[5]{H\left(\left.\begin{array}{ll}{#1}&{#2}\\{#3}&{#4}\end{array}\right|{#5}\right)}
\nc{\web}[5]{\overline{W}\left(\left.\begin{array}{ll}{#1}&{#2}\\{#3}&{#4}\end{array}\right|{#5}\right)}
\nc{\wep}[5]{W'\left(\left.\begin{array}{ll}{#1}&{#2}\\{#3}&{#4}\end{array}\right|{#5}\right)}
\nc{\wes}[5]{W^*\left(\left.\begin{array}{ll}{#1}&{#2}\\{#3}&{#4}\end{array}\right|{#5}\right)}
\nc{\wess}[5]{W^{**}\left(\left.\begin{array}{ll}{#1}&{#2}\\{#3}&{#4}\end{array}\right|{#5}\right)}
\nc{\cet}[7]{C^{#6}_{#7}\left(\left.\begin{array}{ll}{#1}&{#2}\\{#3}&{#4}\end{array}\right|{#5}\right)}
\nc{\bcet}[7]{\bar{C}^{#6}_{#7}\left(\left.\begin{array}{ll}{#1}&{#2}\\{#3}&{#4}\end{array}\right|{#5}\right)}
\nc{\wet}[7]{W^{#6}_{#7}\left(\left.\begin{array}{ll}{#1}&{#2}\\{#3}&{#4}\end{array}\right|{#5}\right)}
\nc{\bwet}[7]{\overline{W}^{#6}_{#7}\left(\left.\begin{array}{ll}{#1}&{#2}\\{#3}&{#4}\end{array}\right|{#5}\right)}
\nc{\wec}[7]{\widetilde{W}^{#6}_{#7}\left(\left.\begin{array}{ll}{#1}&{#2}\\{#3}&{#4}\end{array}\right|{#5}\right)}
\nc{\wgen}[6]{W^{#6}\left(\left.\begin{array}{ll}{#1}&{#2}\\{#3}&{#4}\end{array}\right|{#5}\right)}
\nc{\wgenp}[6]{W^{*{#6}}\left(\left.\begin{array}{ll}{#1}&{#2}\\{#3}&{#4}\end{array}\right|{#5}\right)}
\nc{\wo}[5]{\Omega\left(\left.\begin{array}{ll}{#1}&{#2}\\{#3}&{#4}\end{array}\right|{#5}\right)}
\nc{\wsgen}[8]{{#8}^{#6}_{#7}\left(\left.\begin{array}{ll}{#1}&{#2}\\{#3}&{#4}\end{array}\right|{#5}\right)}
\nc{\qbinom}[2]{{\genfrac{[}{]}{0pt}{}{{#1}}{{#2}}}_{q}}
\nc{\hhg}[4]{\phi\left({{{#1}\,\,\,{#2}}\atop{{#3}}};
                     {#4}\right)}
\nc{\fullhhg}[5]{{_1}\phi_2\left({{{#1}\,\,\,{#2}}\atop{{#3}}};
                     {#4},{#5}\right)}
\nc{\qp}[2]{({#1}\, ; \, {#2})_{\infty}}
\nc{\qpf}[1]{({#1}\, ; \, q^4)_{\infty}}
\nc{\pp}[1]{({#1}\, ; \, p)_{\infty}}
\nc{\qpp}[1]{({#1}\, ; \, p, q^4)_{\infty}}
\begin{document}

\ 
%\vspace{5mm}
\begin{center}
{\LARGE\bf The Vertex-Face Correspondence \\
 and the Elliptic $6j$-symbols}

\vspace{25mm}
{\large Hitoshi Konno}

\vspace{5mm}
{\large\it 
Department of Mathematics, Faculty of Integrated Arts $\&$ Sciences, \\
Hiroshima University, Higashi Hiroshima 739-8521, Japan\footnote{Email: konno@mis.hiroshima-u.ac.jp}}
%\date{
%{\rm \today}
%}

\vspace{30mm}
{\bf Abstract}
\end{center}

\noindent
A new formula connecting the elliptic $6j$-symbols and the fusion of 
the vertex-face intertwining vectors is given. This is based on 
the identification of the $k$ fusion 
intertwining vectors with the change of base matrix elements from 
Sklyanin's standard base to Rosengren's 
natural base in the space of even theta functions of order $2k$. 
The new formula allows us to derive various properties of  the elliptic 
$6j$-symbols, 
such as the addition formula, the biorthogonality property, the fusion formula and the Yang-Baxter relation.
We also discuss a connection with the Sklyanin algebra based on the factorised formula for 
the $L$-operator.

\newpage
\section{Introduction}
The theory of elliptic hypergeometric series has been rapidly developing. 
In \cite{FreTur}, Frenkel and Turaev took the initiative 
introducing a notion of modular hypergeometric functions ${}_{r+1}{\omega}_r$ 
as an elliptic analogue of the very-well-poised balanced basic hypergeometric series 
${}_{r+1}\varphi_r$. This was based on an observation that the fusion of the 
face type Boltzmann weights studied by Date et.al.\cite{DJKMO} can be regarded as 
an elliptic analogue of the $q$-$6j$ symbols.

The theory of biorthogonal rational functions led Spiridonov and Zhedanov
 to a generalization of
the elliptic $6j$-symbols\cite{SpiZhed}. Deeper consideration of the well-poisedness
and very-well-poisedness conditions for general elliptic hypergeometric
series resulted in a change of notations to $ _{r+3}V_{r+2}$ \cite{Spi}, which is
now considered as a proper elliptic analogue of the very-well-poised basic
hypergeometric series $ _{r+1}\varphi_r $.
At the same time, they reached a relevant scheme for dealing with a 
family of biorthogonal rational functions as 
the generalized eigenvalue problem (GEVP) associated with two Jacobi matrices.  
%Accordingly a notion of ``elliptic $6j$-symbols'' was generalised to the one 
% should be expressed by ${}_{12}V_{11}$.  
 
Recently Rosengren  found 
a GEVP which is relevant specific to a generalisation of the  
elliptic $6j$-symbols\cite{Ro1,Ro2}. 
It is deeply related with the representation of the Sklyanin algebra\cite{Sklyanin}
 on the space of  theta functions. 
He found a {\it natural basis} of the space as a set of solutions of the GEVP
and identified a change of base matrix elements between the two natural bases 
depending on a different parameters with 
a generalisation of the elliptic $6j$-symbols.
 He then succeeded to derive an expression of the generalised ``elliptic $6j$-symbols'' in terms of ${}_{12}V_{11}$. 

Such a connection with the Sklyanin algebra reminds us of the 
work by Takebe\cite{Takebe} 
on the diagonalisation of the higher spin eight-vertex 
model by using the algebraic Bethe ansatz. In fact, Rosengren's natural basis 
turned out to be a fusion of the vertex-face intertwining vectors realised as the vectors in the space of theta functions. 
Such vectors are used to map (or gauge transform) 
a higher dimensional representation of the Sklyanin algebra 
associated with the fusion of the  elliptic $R$-matrix to the Bethe ansatz operators of 
the fusion of the eight-vertex 
Solid-on-Solid (SOS) model. This makes the diagonalisation of the model simpler. 
The eigen states (the Bethe states) are then given by the fused intertwining 
vectors. 

The purpose of this paper is to make a direct connection 
between the generalised 
elliptic 
$6j$-symbols and the fusion of the vertex-face intertwining vectors. 
In a context of the solvable 
lattice models, it is standard to take intertwining vectors as  
the elements of the vector space $V$, on which tensor product $V\otimes V$ the 
fusion of the $R$-matrix acts. They map a fusion of the elliptic $R$-matrix 
 to a fusion of the face weight of the SOS model. 
By making a connection between Takebe's  realisation of the intertwining vectors 
and the standard one, we find an 
exact relation between Rosengren's change of base matrix elements 
and the fusion of the standard intertwining vectors. For this purpose, 
we fully use the dual of the intertwining vectors and 
their fusions. In the process, we 
also give a formula connecting 
the elliptic binomial coefficients and the intertwining vectors.

This paper is organised as follows. In the next section, we make a quick review of the fusion of the eight-vertex model and the eight-vertex SOS model as well as their vertex-face correspondence relationship. These follow the results in \cite{DJKMO,KKW}. 
In Sec.3, we make a connection 
between the standard intertwining vectors and Rosengren's vectors of 
natural basis. This leads us to a new formula for the generalised elliptic $6j$-symbols. We then derive various properties of the elliptic $6j$-symbols. We also make a brief 
comment on a connection with the Sklyanin algebra. 
Since we use a slightly different notation from \cite{Ro1,Ro2}, in the Appendix, 
we give a derivation of the elliptic binomial theorem, 
the elliptic Jackson's summation formula and the expression of 
the elliptic $6j$-symbols in terms of ${}_{12}V_{11}$ following the idea of Rosengren. The other idea of proof can be found in \cite{KN,SpiZhed2}.

\subsection{Notations}

Let ${p}=e^{-\frac{\pi K'}{K}}$, $q=-e^{-\frac{\pi \lambda}{2K}}$ and $\zeta=e^{-\frac{\pi \la u}{2K}}$.
We introduce $x$, $\tau$ and $r$ by $x=-q$, 
$\tau=\frac{2iK}{K'}$ and $r=\frac{K'}{\la}$. Then $p=e^{-\frac{2\pi i}{ \tau}}=x^{2r}$. The parameter $r$ plays a role of restriction height in the restricted SOS models. Through this paper, we assume ${\rm Im}\tau >0$.
Let $\tp=e^{2\pi i \tau}$.
%=e^{-\pi \frac{I'}{I}}$, where $I=\frac{K'}{2},\ I'=2K$. 
We use the theta functions 
\be
&&\vartheta_1(u|\tau)=2\tilde{p}^{1/8}(\tilde{p};\tilde{p})_\infty\sin\pi u
\prod_{n=1}^\infty(1-2\tilde{p}^n\cos2\pi u+\tilde{p}^{2n}),\\
&&\vartheta_0(u|\tau)=-ie^{\pi i(u+\tau/4)}\vartheta_1\left(u+\frac{\tau}{2}\Big|\tau\right),\\
&&\vartheta_2(u|\tau)=\vartheta_1\left(u+\frac{1}{2}\Big|\tau\right),\\
&&\vartheta_3(u|\tau)=e^{\pi i(u+\tau/4)}\vartheta_1\left(u+\frac{\tau+1}{2}\Big|\tau\right).
\en
We also use the symbol $[u]$ defined by 
\be
&&[u]=x^{\frac{u^2}{r}-u}\Theta_{x^{2r}}(x^{2u})=C\vt{1}{u}{\tau}, \quad C=x^{-\frac{r}{4}}e^{-\frac{\pi i}{4}}\tau^{\frac{1}{2}}.
\en
The elliptic shifted factorials are defined by 
\be
&&[u]_n=\prod_{j=0}^{n-1}[u+j]
\en
with the convention
\be
&&[u_1,u_2,\cdots,u_k]_n=\prod_{i=1}^k[u_i]_n.
\en

\section{Fusion and the Vertex-Face Correspondence}
According to \cite{DJKMO,KKW}, we here make a brief review of the fusion of the 
Boltzmann weight associated with 
the eight-vertex model and the eight-vertex solid-on-solid (SOS) model as well as 
 their vertex-face correspondence relationship. 

\subsection{Fusion of Baxter's elliptic $R$-matrix}
The eight-vertex model is a two-dimensional square lattice model. 
For each vertex, we associate a dynamical variable (a spin) 
$\vep_j\in \{+, -\}$ with each edge $j$ (Figure 1). 
\begin{center}
\setlength{\unitlength}{0.008mm}%
\begingroup\makeatletter\ifx\SetFigFont\undefined%
\gdef\SetFigFont#1#2#3#4#5{%
  \reset@font\fontsize{#1}{#2pt}%
  \fontfamily{#3}\fontseries{#4}\fontshape{#5}%
  \selectfont}%
\fi\endgroup%
\begin{picture}(4212,3648)(601,-3340)
\put(-551,-3286){\makebox(0,0)[lb]{\smash{\SetFigFont{12}{14.4}{\rmdefault}{\mddefault}{\updefault}{\color[rgb]{0,0,0}{\footnotesize Figure 1: The vertex model weight} }%
}}}
\thinlines
{\color[rgb]{0,0,0}\put(4801,-961){\vector(-1, 0){2400}}
}%
\put(4526,-811){\makebox(0,0)[lb]{\smash{\SetFigFont{12}{14.4}{\rmdefault}{\mddefault}{\updefault}{\color[rgb]{0,0,0}$\vep_2$}%
}}}
\put(3751,164){\makebox(0,0)[lb]{\smash{\SetFigFont{12}{14.4}{\rmdefault}{\mddefault}{\updefault}{\color[rgb]{0,0,0}$\vep_1$}%
}}}
\put(3751,-2000){\makebox(0,0)[lb]{\smash{\SetFigFont{12}{14.4}{\rmdefault}{\mddefault}{\updefault}{\color[rgb]{0,0,0}$\vep_1'$}%
}}}
\put(2501,-811){\makebox(0,0)[lb]{\smash{\SetFigFont{12}{14.4}{\rmdefault}{\mddefault}{\updefault}{\color[rgb]{0,0,0}$\vep_2'$}%
}}}
\put(2001,-961){\makebox(0,0)[lb]{\smash{\SetFigFont{12}{14.4}{\rmdefault}{\mddefault}{\updefault}{\color[rgb]{0,0,0}$v$}%
}}}
\put(3601,-2461){\makebox(0,0)[lb]{\smash{\SetFigFont{12}{14.4}{\rmdefault}{\mddefault}{\updefault}{\color[rgb]{0,0,0}$u$}%
}}}
\put(-1500,-1036){\makebox(0,0)[lb]{\smash{\SetFigFont{12}{14.4}{\rmdefault}{\mddefault}{\updefault}{\color[rgb]{0,0,0}$R(u-v)^{\vep_1\vep_2}_{\vep_1'\vep_2'}$}%
}}}
\put(1426,-1036){\makebox(0,0)[lb]{\smash{\SetFigFont{12}{14.4}{\rmdefault}{\mddefault}{\updefault}{\color[rgb]{0,0,0}=}%
}}}
{\color[rgb]{0,0,0}\put(3601,239){\vector( 0,-1){2400}}
}%
\end{picture}{}\qquad\qquad\qquad{}
\end{center}
Let $V=\C v_{\vep_1}\oplus \C v_{\vep_2}$. 
We regard $R(u)\in \End(V\otimes V)$ by 
\be
&&R(u)v_{\vep_1}\otimes v_{\vep_2}=\sum_{\vep_1',\vep_2'}
R(u)_{\vep_1'\vep_2'}^{\vep_1\vep_2}v_{\vep_1'}\otimes v_{\vep_2'}.
\en
We have only the eight possible configurations, $R_{++}^{++}, 
R_{--}^{--}, R_{++}^{--}, R_{--}^{++}, R_{+-}^{+-}$, $R_{-+}^{-+}$, $R_{+-}^{-+}$, $R_{-+}^{+-}$. To each of them, we assign   
 the corresponding matrix element of Baxter's elliptic $R$-matrix as the Boltzmann weight\cite{Baxter}.  
\bea
&&{R}(u)=
{R}_0(u)\left(\matrix{a(u)&&&d(u)\cr
&b(u)&c(u)&\cr
&c(u)&b(u)&\cr 
d(u)&&&a(u)\cr}\right),\lb{ellR}
\ena
where
\be
R_0(u)&=&z^{-\frac{r-1}{2r}}\frac{(px^2z;x^4,p)_\infty(x^2z;x^4,p)_\infty
(p/z;x^4,p)_\infty(x^4/z;x^4,p)_\infty}{(px^2/z;x^4,p)_\infty
(x^2/z;x^4,p)_\infty(pz;x^4,p)_\infty(x^4z;x^4,p)_\infty},\lb{evR}\\
a(u)&=&\frac{\vtf{2}{\frac{1}{2r}}{\frac{\tau}{2}}\vtf{2}{\frac{  u}{2r}}{\frac{\tau}{2}}}{\vtf{2}{0}{\frac{\tau}{2}}\vtf{2}{\frac{1+u}{2r}}{\frac{\tau}{2}}},\qquad
b(u)=\frac{\vtf{2}{\frac{1 }{2r}}{\frac{\tau}{2}}\vtf{1}{\frac{  u}{2r}}{\frac{\tau}{2}}}
{\vtf{2}{0}{\frac{\tau}{2}}\vtf{1}{\frac{1+u}{2r}}{\frac{\tau}{2}}},\nn\\
c(u)&=&\frac{\vtf{1}{\frac{ 1}{2r}}{\frac{\tau}{2}}\vtf{2}{\frac{  u}{2r}}{\frac{\tau}{2}}}
{\vtf{2}{0}{\frac{\tau}{2}}\vtf{1}{\frac{1+u}{2r}}{\frac{\tau}{2}}},\qquad
d(u)=-\frac{\vtf{1}{\frac{ 1}{2r}}{\frac{\tau}{2}}\vtf{1}{\frac{  u}{2r}}{\frac{\tau}{2}}}
{\vtf{2}{0}{\frac{\tau}{2}}\vtf{2}{\frac{1+u}{2r}}{\frac{\tau}{2}}}\nn
\en
with $z=\zeta^2=x^{2u}$.
  Note that the $R$-matrix $R(u)$ is the two dimensional representation of the 
universal $R$-matrix of the vertex type elliptic quantum group $\Aqp(\slth)$ \cite{JKOSI}. 
The $R$-matrix \eqref{ellR} satisfies the Yang-Baxter equation (YBE), the unitarity relation, the crossing symmetry relation and the initial condition given as follows. 
\bea
&&R_{12}(u-v)R_{13}(u)R_{23}(v)=R_{23}(v)R_{13}(u)R_{12}(u-v),\lb{ybeR}\\
&&R(u)PR(u)P=\id,\\
&&R(-u-1)=(\sigma^{y }\otimes 1)^{-1}\ (PR(u)P)^{t_1} \  \sigma^y\otimes 1,\lb{crossR}\\
&&R(0)=P, \qquad \lim_{u\to -1}R(u)=P-\id.
\ena
Here ${^{t_1}}$ denotes the transposition with respect to the first vector space in the tensor product and $P(v_{\vep_1}\otimes v_{\vep_2})=v_{\vep_2}\otimes v_{\vep_1}$.

Fusion of $R(u)$ was considered systematically in \cite{DJKMO} ( see also \cite{Konno04}  for the $2\times 2$ fusion case). 
Let $V_j, V_{\baj}$ be copies of V. 
Let us define the operator $\Pi_{1 \cdots k}$ by
\be
&&\Pi_{1\cdots k}=\frac{1}{k!}(P_{1k}+\cdots+P_{k-1k}+I)\cdots(P_{13}+P_{23}+I) (P_{12}+I).
\en
This yields the projection on the space $V^{(k)}$ of the symmetric tensors in $V^{\otimes k}$.
We define  
\be
&&R^{(k,1)}_{1\cdots k, \baj}(u)=\Pi_{1\cdots k}R_{1\baj}(u+k-1)\cdots R_{k-1\baj}(u+1)R_{k\baj}(u)\ \in {\rm End}(V^{(k)} \otimes V_{\baj}).
\en
The $R$-matrix $R^{(k,1)}_{1\cdots k, \baj}(u)$ then satisfies 
\be
&&R^{(k,1)}_{1\cdots k, \baj}(u)\Pi_{1\cdots k}=R^{(k,1)}_{1\cdots k, \baj}(u).
\en
By fusing $R^{(k,1)}_{1\cdots k, \baj}(u)$, $l$ times, we then define 
the  $k\times l$ fusion $R$-matrix as follows.
\bea
&&R^{(k,l)}(u)=\Pi_{\bao\cdots \bal}R^{(k,1)}_{1\cdots k, \bal}(u)R^{(k,1)}_{1\cdots k, \overline{l-1}}(u-1)\cdots R^{(k,1)}_{1\cdots k, \bao}(u-l+1). \lb{kkfusionR}
\ena
This is an operator in ${\rm End}(V^{(k)} \otimes V^{(l)})$. It
satisfies 
\bea
&&R^{(k,l)}(u)=R^{(k,l)}(u)\Pi_{1\cdots k}=R^{(k,l)}(u)\Pi_{\bao\cdots \bal}.\lb{projkk}
\ena
Using the YBE \eqref{ybeR} repeatedly and \eqref{projkk}, we verify that $R^{(k,l)}(u)$ satisfies the YBE 
on $V^{(k)}\otimes V^{(l)} \otimes V^{(m)}$. 
\be
&&R^{(k,l)}(u-v)R^{(k,m)}(u)R^{(l,m)}(v)=R^{(l,m)}(v)R^{(k,m)}(u)R^{(k,l)}(u-v).
\en

\subsection{Fusion of the eight-vertex SOS model}

The eight-vertex SOS model is also a two-dimensional square lattice model. The dynamical variables $a_j$ are called local heights.
They take values in $\Z$. For each face, we associate a local height $a_j$ with each vertex $j$ (Figure 2). 
\be
{}\qquad \qquad \setlength{\unitlength}{0.01mm}%
\begingroup\makeatletter\ifx\SetFigFont\undefined%
\gdef\SetFigFont#1#2#3#4#5{%
  \reset@font\fontsize{#1}{#2pt}%
  \fontfamily{#3}\fontseries{#4}\fontshape{#5}%
  \selectfont}%
\fi\endgroup%
\begin{picture}(4500,3264)(526,-3040)
\put(-501,-2986){\makebox(0,0)[lb]{\smash{\SetFigFont{12}{14.4}{\rmdefault}{\mddefault}
{\updefault}{\color[rgb]{0,0,0} {\footnotesize Figure 2: The SOS model face weight } }%
}}}
\thinlines
{\color[rgb]{0,0,0}\put(3001,-61){\line( 1, 0){1800}}
\put(4801,-61){\line( 0,-1){1800}}
\put(4801,-1861){\line(-1, 0){1800}}
}%
{\color[rgb]{0,0,0}\put(3501,-61){\line(-1,-1){500}}
}%
\put(3726,-1050){\makebox(0,0)[lb]{\smash{\SetFigFont{12}{14.4}{\rmdefault}{\mddefault}{\updefault}{\color[rgb]{0,0,0}$u$}%
}}}
\put(5026,-2011){\makebox(0,0)[lb]{\smash{\SetFigFont{12}{14.4}{\rmdefault}{\mddefault}{\updefault}{\color[rgb]{0,0,0}$a_3$}%
}}}
\put(4951, 89){\makebox(0,0)[lb]{\smash{\SetFigFont{12}{14.4}{\rmdefault}{\mddefault}{\updefault}{\color[rgb]{0,0,0}$a_2$}%
}}}
\put(2701, 89){\makebox(0,0)[lb]{\smash{\SetFigFont{12}{14.4}{\rmdefault}{\mddefault}{\updefault}{\color[rgb]{0,0,0}$a_1$}%
}}}
\put(2701,-2086){\makebox(0,0)[lb]{\smash{\SetFigFont{12}{14.4}{\rmdefault}{\mddefault}{\updefault}{\color[rgb]{0,0,0}$a_4$}%
}}}
\put(-1026,-961){\makebox(0,0)[lb]{\smash{\SetFigFont{12}{14.4}{\rmdefault}{\mddefault}{\updefault}{\color[rgb]{0,0,0}$W\BW{a_1}{a_2}{a_4}{a_3}{u}$}%
}}}
\put(2101,-961){\makebox(0,0)[lb]{\smash{\SetFigFont{12}{14.4}{\rmdefault}{\mddefault}{\updefault}{\color[rgb]{0,0,0}=}%
}}}
{\color[rgb]{0,0,0}\put(3001,-61){\line( 0,-1){1800}}
}%
%\put(3001,-450){\line(1,1){400}}
\lb{SOS}\end{picture}
\en
We allow only the configurations satisfying the so-called the admissibility condition $|a_j-a_k|=1$ for 
any two adjacent local heights $a_j$ and $a_k$. Then we have only the six possible 
configurations for each face. 
 We assign the following     
Boltzmann weight (or face weight) $W\BW{a_1}{a_2}{a_4}{a_3}{u}$ to each configuration.  
  \begin{eqnarray}
&&{W}\left(\left.
\begin{array}{cc}
a&a\pm 1\\
a\pm1&a\pm2
\end{array}\right|u\right)={R}_0(u),\nn\\
&&{W}\left(\left.
\begin{array}{cc}
a&a\pm 1\\
a\pm1&a
\end{array}\right|u\right)={R}_0(u)\frac{[a\mp u][1]}{[ a][ 1+u]},\lb{face}\\
&&{W}\left(\left.
\begin{array}{cc}
a&a\pm 1\\
a\mp 1&a
\end{array}\right|u\right)=
{R}_0(u)
\frac{[ a\pm 1][u] }{[ a][1+u]}.\nn
\end{eqnarray}
Note that these weights can be understood as the matrix elements of 
the dynamical $R$-matrix 
$R(u,a)$ obtained as the two dimensional 
representation of the 
universal dynamical $R$-matrix of the face type elliptic quantum group 
$\Bqla(\slth)$ \cite{JKOSI,JKOSII}. ${\cal B}_{q,\la}(\slth)$ is a 
central extension of Felder's elliptic quantum group\cite{Felder}. 
%The elliptic algebra $U_{q,p}(\slth)$ gives an elliptic analogue of 
%the Drinfeld realisation of $U_q(\slth)$ for $\Bqla(\slth)$\cite{Konno,JKOSII}.

The face weights enjoy the following face-type YBE,
unitarity,
and 
crossing symmetry relations.
\bea
&&\sum_{g}W\BW{a}{b}{f}{g}{u}W\BW{f}{g}{e}{d}{v}W\BW{b}{c}{g}{d}{u-v}\nn\\
&&\qquad =\sum_{g}W\BW{a}{g}{f}{e}{u-v}W\BW{a}{b}{g}{c}{v}W\BW{g}{c}{e}{d}{u},\lb{fyb}\\
&&\sum_{e}W\BW{a}{b}{e}{c}{u}W\BW{a}{e}{d}{c}{-u}=1,\lb{fcross}\\
&&W\BW{a}{b}{d}{c}{u}=(-)^{\frac{a+d-b-c}{2}}\frac{[b]}{[a]}W\BW{d}{a}{c}{b}{-u-1}.\lb{crossW}
\ena
In addition, the following initial conditions are the key 
for the fusion of the SOS weights.
\be
&&W\BW{a}{b}{d}{c}{0}=\delta_{b,d},\\
&&W\BW{a}{b}{d}{c}{-1}=0\qquad {\rm if}\ |a-c|=2,\\
&&W\BW{a}{a\pm 1}{a\pm 1}{a}{-1}=-W\BW{a}{a\pm 1}{a\mp 1}{a}{-1}.
\en

The  $k\times l$ fusion of the face weight is obtained following the two steps\cite{DJKMO}. 
Firstly we define
\be
&&W^{(k,1)}\BW{a}{b}{d}{c}{u}\nn\\
&&=\sum_{d_1,..,d_{k-1}}  \!\!\!\! W\BW{a}{a_1}{d}{d_1}{u+k-1}W\BW{a_1}{b}{d_1}{c}{u+k-2}\cdots
W\BW{a_{k-1}}{b}{d_{k-1}}{c}{u}.\nn\\
\en
Then the RHS is independent of the choice of $a_1,..,a_{k-1}$ provided
$|a-a_1|=|a_1-a_2|=\cdots=|a_{k-1}-b|=1$. Secondly we define
\bea
&&
%\hspace*{-20mm}
W^{(k,l)}\BW{a}{b}{d}{c}{u}\nn\\
&&
%\hspace*{-30mm}
=\sum_{a_1,..,a_{l-1}}W^{(k,1)}\BW{a}{b}{a_1}{b_1}{u-l+1}
\!\! W^{(k,1)}\BW{a}{b}{a_1}{b_1}{u-l+2}
\cdots W^{(k,1)}\BW{a_{l-1}}{b_{l-1}}{d}{c}{u}. \nn\\
\lb{kkfusionW}
\ena
Then the RHS is independent of the choice of $b_1,..,b_{l-1}$ provided $|b-b_1|=|b_1-b_2|=\cdots=|b_{l-1}-c|=1$.
In $W^{(k,l)}$, the admissible condition for the dynamical variables is extended to $a-b\ \in \{-k, -k+2, .., k\}$ for any two horizontally adjacent local heights $a, b$, while 
$a-d\ \in \{-l, -l+2, .., l\}$ for any two vertically adjacent local heights $a, d$. 
The $k\times l$ fusion face weight $W^{(k,l)}$ satisfies the face type YBE.
\bea
&&\sum_{g}W^{(k,l)}\BW{a}{b}{f}{g}{u}W^{(k,m)}\BW{f}{g}{e}{d}{v}W^{(m,l)}\BW{b}{c}{g}{d}{u-v}\nn\\
&&\qquad =\sum_{g}W^{(m,l)}\BW{a}{g}{f}{e}{u-v}W^{(k,m)}\BW{a}{b}{g}{c}{v}W^{(k,l)}\BW{g}{c}{e}{d}{u}.\lb{fyb}
\ena

\subsection{The vertex-face correspondence }
The vertex-face correspondence is a relationship between Baxter's $R$-matrix and the SOS face weight $W\BW{a_1}{a_2}{a_4}{a_3}{u}$. 
Let us consider the following vectors in $V$   
\bea
&&\psi(u)^a_b=\psi_+(u)^a_b\ v_+ + \psi_-(u)^a_b\ v_-,\qquad \lb{intertwinvec}\\
&&{\psi}_+(u)^a_b=\vtf{0}{\frac{(a-b)u+a}{2r}}{\frac{\tau}{2}},\qquad \psi_-(u)^a_b=\vtf{3}{\frac{(a-b)u+a}{2r}}{\frac{\tau}{2}} \nn
\ena 
with $|a-b|=1$. Baxter showed that the following identity holds\cite{Baxter}.
\bea
&&\sum_{\vep_1',\vep_2'}
R(u-v)_{\vep_1 \vep_2}^{\vep_1' \vep_2'}\ 
\psi_{\vep_1'}(u)_{b}^{a}
\psi_{\vep_2'}(v)_{c}^{b}
=\sum_{b' \in {\mathbb{Z}}}
\psi_{\vep_2}(v)_{b'}^{a}
\psi_{\vep_1}(u)_{c}^{b'}
W\left(\left.
\begin{array}{cc}
a&b\\
b'&c
\end{array}\right|u-v\right).\nn\\
&&\lb{vertexface}
\ena
We call $\psi(u)^a_b$ the (vertex-face) intertwining vectors. 

\vspace{5mm}
\begin{center}
\setlength{\unitlength}{0.01mm}%
\begingroup\makeatletter\ifx\SetFigFont\undefined%
\gdef\SetFigFont#1#2#3#4#5{%
  \reset@font\fontsize{#1}{#2pt}%
  \fontfamily{#3}\fontseries{#4}\fontshape{#5}%
  \selectfont}%
\fi\endgroup%
\begin{picture}(9312,2889)(301,-2815)
\put(7651,-2761){\makebox(0,0)[lb]{\smash{\SetFigFont{12}{14.4}{\rmdefault}{\mddefault}{\updefault}{\color[rgb]{0,0,0}$(b)$}%
}}}
\thinlines
{\color[rgb]{0,0,0}\put(3001,-361){\vector( 0,-1){1200}}
}%
{\color[rgb]{0,0,0}\put(8401,-1561){\line( 1, 0){1200}}
}%
{\color[rgb]{0,0,0}\put(9001,-361){\vector( 0,-1){600}}
}%
{\color[rgb]{0,0,0}\put(9001,-961){\line( 0,-1){600}}
}%
\put(2101,-61){\makebox(0,0)[lb]{\smash{\SetFigFont{12}{14.4}{\rmdefault}{\mddefault}{\updefault}{\color[rgb]{0,0,0}$a$}%
}}}
\put(3451,-61){\makebox(0,0)[lb]{\smash{\SetFigFont{12}{14.4}{\rmdefault}{\mddefault}{\updefault}{\color[rgb]{0,0,0}$b$}%
}}}
\put(3226,-1561){\makebox(0,0)[lb]{\smash{\SetFigFont{12}{14.4}{\rmdefault}{\mddefault}{\updefault}{\color[rgb]{0,0,0}$\vep$}%
}}}
\put(2901,-1861){\makebox(0,0)[lb]{\smash{\SetFigFont{12}{14.4}{\rmdefault}{\mddefault}{\updefault}{\color[rgb]{0,0,0}$u$}%
}}}
\put(201,-961){\makebox(0,0)[lb]{\smash{\SetFigFont{12}{14.4}{\rmdefault}{\mddefault}{\updefault}{\color[rgb]{0,0,0}$\psi_\vep(u)^a_b$}%
}}}
\put(1501,-961){\makebox(0,0)[lb]{\smash{\SetFigFont{12}{14.4}{\rmdefault}{\mddefault}{\updefault}{\color[rgb]{0,0,0}$=$}%
}}}
\put(1801,-2761){\makebox(0,0)[lb]{\smash{\SetFigFont{12}{14.4}{\rmdefault}{\mddefault}{\updefault}{\color[rgb]{0,0,0}$(a)$}%
}}}
\put(7501,-961){\makebox(0,0)[lb]{\smash{\SetFigFont{12}{14.4}{\rmdefault}{\mddefault}{\updefault}{\color[rgb]{0,0,0}$=$}%
}}}
\put(8901,-61){\makebox(0,0)[lb]{\smash{\SetFigFont{12}{14.4}{\rmdefault}{\mddefault}{\updefault}{\color[rgb]{0,0,0}$u$}%
}}}
\put(9151,-661){\makebox(0,0)[lb]{\smash{\SetFigFont{12}{14.4}{\rmdefault}{\mddefault}{\updefault}{\color[rgb]{0,0,0}$\vep$}%
}}}
\put(8326,-1861){\makebox(0,0)[lb]{\smash{\SetFigFont{12}{14.4}{\rmdefault}{\mddefault}{\updefault}{\color[rgb]{0,0,0}$a$}%
}}}
\put(9601,-1861){\makebox(0,0)[lb]{\smash{\SetFigFont{12}{14.4}{\rmdefault}{\mddefault}{\updefault}{\color[rgb]{0,0,0}$b$}%
}}}
\put(6001,-961){\makebox(0,0)[lb]{\smash{\SetFigFont{12}{14.4}{\rmdefault}{\mddefault}{\updefault}{\color[rgb]{0,0,0}$\psi^*_\vep(u)_a^b$}%
}}}
{\color[rgb]{0,0,0}\put(2401,-361){\line( 1, 0){1200}}
}%
\end{picture}

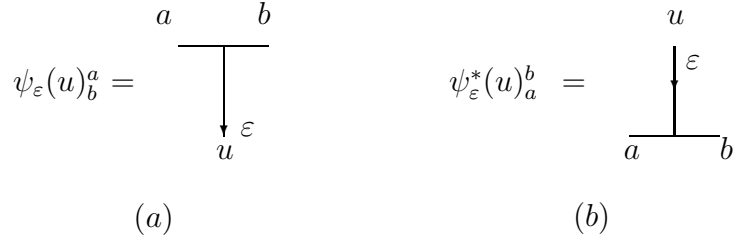
{\footnotesize Figure 3: \lb{ITV} $(a)$ The intertwining vector ; $(b)$ the dual intertwining vector}
\end{center}

\vspace{5mm}
\begin{center}
\setlength{\unitlength}{0.012mm}%
\begingroup\makeatletter\ifx\SetFigFont\undefined%
\gdef\SetFigFont#1#2#3#4#5{%
  \reset@font\fontsize{#1}{#2pt}%
  \fontfamily{#3}\fontseries{#4}\fontshape{#5}%
  \selectfont}%
\fi\endgroup%
\begin{picture}(11850,3369)(451,-3190)
\put(9001,-3136){\makebox(0,0)[lb]{\smash{\SetFigFont{12}{14.4}{\rmdefault}{\mddefault}{\updefault}{\color[rgb]{0,0,0}$(b)$}%
}}}
\thinlines
{\color[rgb]{0,0,0}\put(1801,-361){\vector( 0,-1){1500}}
}%
{\color[rgb]{0,0,0}\put(2401,-961){\vector(-1, 0){1500}}
}%
{\color[rgb]{0,0,0}\put(8401,-1561){\line(-1, 0){1200}}
\put(7201,-1561){\line( 0, 1){1200}}
}%
{\color[rgb]{0,0,0}\put(10501,-361){\line( 0,-1){1200}}
\put(10501,-1561){\line( 1, 0){1200}}
\put(11701,-1561){\line( 0, 1){1200}}
\put(11701,-361){\line(-1, 0){1200}}
}%
{\color[rgb]{0,0,0}\put(7801,-136){\vector( 0,-1){1500}}
}%
{\color[rgb]{0,0,0}\put(8701,-961){\vector(-1, 0){1500}}
}%
{\color[rgb]{0,0,0}\put(11101,-61){\vector( 0,-1){300}}
}%
{\color[rgb]{0,0,0}\put(4501,-361){\line( 0,-1){1200}}
\put(4501,-1561){\line( 1, 0){1200}}
\put(5701,-1561){\line( 0, 1){1200}}
\put(5701,-361){\line(-1, 0){1200}}
}%
{\color[rgb]{0,0,0}\put(4501,-961){\vector(-1, 0){300}}
}%
{\color[rgb]{0,0,0}\put(5101,-1561){\vector( 0,-1){300}}
}%
{\color[rgb]{0,0,0}\put(12001,-961){\vector(-1, 0){300}}
}%
{\color[rgb]{0,0,0}\put(4801,-361){\line(-1,-1){300}}
}%
{\color[rgb]{0,0,0}\put(10801,-361){\line(-1,-1){300}}
}%
\put(3301,-961){\makebox(0,0)[lb]{\smash{\SetFigFont{12}{14.4}{\rmdefault}{\mddefault}{\updefault}{\color[rgb]{0,0,0}$=$ }%
}}}
\put(1651,-2311){\makebox(0,0)[lb]{\smash{\SetFigFont{12}{14.4}{\rmdefault}{\mddefault}{\updefault}{\color[rgb]{0,0,0}$u$}%
}}}
\put(1951,-1861){\makebox(0,0)[lb]{\smash{\SetFigFont{12}{14.4}{\rmdefault}{\mddefault}{\updefault}{\color[rgb]{0,0,0}$\vep_1$}%
}}}
\put(901,-1261){\makebox(0,0)[lb]{\smash{\SetFigFont{12}{14.4}{\rmdefault}{\mddefault}{\updefault}{\color[rgb]{0,0,0}$\vep_2$}%
}}}
\put(451,-961){\makebox(0,0)[lb]{\smash{\SetFigFont{12}{14.4}{\rmdefault}{\mddefault}{\updefault}{\color[rgb]{0,0,0}$v$}%
}}}
\put(901,-211){\makebox(0,0)[lb]{\smash{\SetFigFont{12}{14.4}{\rmdefault}{\mddefault}{\updefault}{\color[rgb]{0,0,0}$a$}%
}}}
\put(2251,-211){\makebox(0,0)[lb]{\smash{\SetFigFont{12}{14.4}{\rmdefault}{\mddefault}{\updefault}{\color[rgb]{0,0,0}$b$}%
}}}
\put(2551,-1561){\makebox(0,0)[lb]{\smash{\SetFigFont{12}{14.4}{\rmdefault}{\mddefault}{\updefault}{\color[rgb]{0,0,0}$c$}%
}}}
\put(4351,-211){\makebox(0,0)[lb]{\smash{\SetFigFont{12}{14.4}{\rmdefault}{\mddefault}{\updefault}{\color[rgb]{0,0,0}$a$}%
}}}
\put(5551,-211){\makebox(0,0)[lb]{\smash{\SetFigFont{12}{14.4}{\rmdefault}{\mddefault}{\updefault}{\color[rgb]{0,0,0}$b$}%
}}}
\put(5851,-1561){\makebox(0,0)[lb]{\smash{\SetFigFont{12}{14.4}{\rmdefault}{\mddefault}{\updefault}{\color[rgb]{0,0,0}$c$}%
}}}
\put(5026,-2311){\makebox(0,0)[lb]{\smash{\SetFigFont{12}{14.4}{\rmdefault}{\mddefault}{\updefault}{\color[rgb]{0,0,0}$u$}%
}}}
\put(3751,-961){\makebox(0,0)[lb]{\smash{\SetFigFont{12}{14.4}{\rmdefault}{\mddefault}{\updefault}{\color[rgb]{0,0,0}$v$}%
}}}
\put(5251,-1936){\makebox(0,0)[lb]{\smash{\SetFigFont{12}{14.4}{\rmdefault}{\mddefault}{\updefault}{\color[rgb]{0,0,0}$\vep_1$}%
}}}
\put(3901,-1261){\makebox(0,0)[lb]{\smash{\SetFigFont{12}{14.4}{\rmdefault}{\mddefault}{\updefault}{\color[rgb]{0,0,0}$\vep_2$}%
}}}
\put(8326,-1936){\makebox(0,0)[lb]{\smash{\SetFigFont{12}{14.4}{\rmdefault}{\mddefault}{\updefault}{\color[rgb]{0,0,0}$a$}%
}}}
\put(6901,-1936){\makebox(0,0)[lb]{\smash{\SetFigFont{12}{14.4}{\rmdefault}{\mddefault}{\updefault}{\color[rgb]{0,0,0}$b$}%
}}}
\put(6901,-211){\makebox(0,0)[lb]{\smash{\SetFigFont{12}{14.4}{\rmdefault}{\mddefault}{\updefault}{\color[rgb]{0,0,0}$c$}%
}}}
\put(7951,-286){\makebox(0,0)[lb]{\smash{\SetFigFont{12}{14.4}{\rmdefault}{\mddefault}{\updefault}{\color[rgb]{0,0,0}$\vep_1$}%
}}}
\put(8251,-1186){\makebox(0,0)[lb]{\smash{\SetFigFont{12}{14.4}{\rmdefault}{\mddefault}{\updefault}{\color[rgb]{0,0,0}$\vep_2$}%
}}}
\put(7726, 89){\makebox(0,0)[lb]{\smash{\SetFigFont{12}{14.4}{\rmdefault}{\mddefault}{\updefault}{\color[rgb]{0,0,0}$u$}%
}}}
\put(8926,-1036){\makebox(0,0)[lb]{\smash{\SetFigFont{12}{14.4}{\rmdefault}{\mddefault}{\updefault}{\color[rgb]{0,0,0}$v$}%
}}}
\put(12301,-1036){\makebox(0,0)[lb]{\smash{\SetFigFont{12}{14.4}{\rmdefault}{\mddefault}{\updefault}{\color[rgb]{0,0,0}$v$}%
}}}
\put(11026, 89){\makebox(0,0)[lb]{\smash{\SetFigFont{12}{14.4}{\rmdefault}{\mddefault}{\updefault}{\color[rgb]{0,0,0}$u$}%
}}}
\put(11326,-136){\makebox(0,0)[lb]{\smash{\SetFigFont{12}{14.4}{\rmdefault}{\mddefault}{\updefault}{\color[rgb]{0,0,0}$\vep_1$}%
}}}
\put(11851,-1186){\makebox(0,0)[lb]{\smash{\SetFigFont{12}{14.4}{\rmdefault}{\mddefault}{\updefault}{\color[rgb]{0,0,0}$\vep_2$}%
}}}
\put(10201,-1936){\makebox(0,0)[lb]{\smash{\SetFigFont{12}{14.4}{\rmdefault}{\mddefault}{\updefault}{\color[rgb]{0,0,0}$b$}%
}}}
\put(10201,-136){\makebox(0,0)[lb]{\smash{\SetFigFont{12}{14.4}{\rmdefault}{\mddefault}{\updefault}{\color[rgb]{0,0,0}$c$}%
}}}
\put(11626,-1936){\makebox(0,0)[lb]{\smash{\SetFigFont{12}{14.4}{\rmdefault}{\mddefault}{\updefault}{\color[rgb]{0,0,0}$a$}%
}}}
\put(9601,-1036){\makebox(0,0)[lb]{\smash{\SetFigFont{12}{14.4}{\rmdefault}{\mddefault}{\updefault}{\color[rgb]{0,0,0}$=$}%
}}}
\put(3001,-3136){\makebox(0,0)[lb]{\smash{\SetFigFont{12}{14.4}{\rmdefault}{\mddefault}{\updefault}{\color[rgb]{0,0,0}$(a)$}%
}}}
{\color[rgb]{0,0,0}\put(1201,-361){\line( 1, 0){1200}}
\put(2401,-361){\line( 0,-1){1200}}
}%
\end{picture}

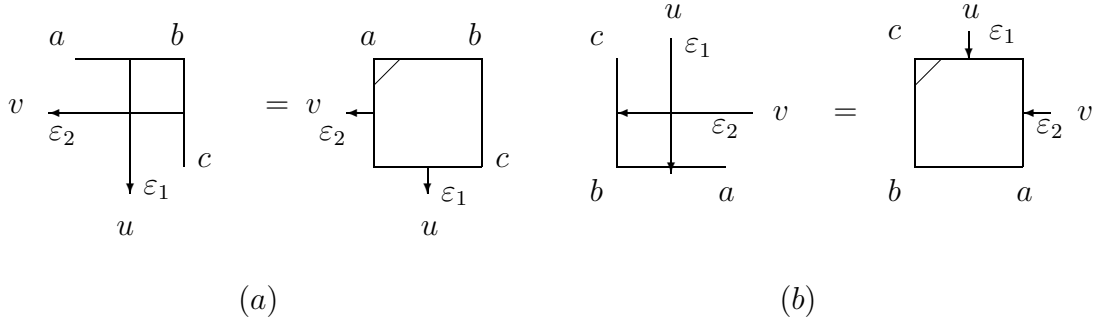
{\footnotesize Figure 4\lb{VFC} The vertex-face correspondence: $(a)$ via the intertwining vector ; $(b)$ via the dual intertwining vector}
\end{center}

In addition, applying the crossing symmetry properties of $R$ \eqref{crossR} and $W$ \eqref{crossW} to \eqref{vertexface}, 
we obtain the following relation. 
\bea
&&\sum_{\vep_1',\vep_2'}
R(u-v)^{\vep_1 \vep_2}_{\vep_1' \vep_2'}\ 
\psi^*_{\vep_1'}(u)_{b}^{a}
\psi^*_{\vep_2'}(v)_{c}^{b}
=\sum_{s \in {\mathbb{Z}}}
\psi^*_{\vep_2}(v)_{b'}^{a}
\psi^*_{\vep_1}(u)_{c}^{b'}
W\left(\left.
\begin{array}{cc}
c&b'\\
b&a
\end{array}\right|u-v\right),\nn\\
&&\lb{vertexfacedual} 
\ena
where we set 
%define $\psi^*_\vep(u)^a_b$  by 
\bea
&&\psi^*_\vep(u)^a_b=-\vep\frac{a-b}{2[b][u]}C^2\ \psi_{-\vep}(u-1)^a_b \lb{dualintvec}
\ena
with $|a-b|=1$. Defining $\psi^*(u)^a_b \in V^*$ by
\be
&&\psi^*(u)^a_b\ v_{\vep}= \psi^*_\vep(u)^a_b, \quad v_{\vep}\in V, 
\en
we call $\psi^*(u)^a_b$ the dual intertwining vectors. In fact, 
by a direct calculation, one can verify that the following
 inversion relations hold.
\bea
&&\sum_{\vep=\pm}\psi_\vep^*(u)^a_b\psi_\vep(u)^b_c=\delta_{a,c},\lb{inversiona}\\
&&\sum_{a=b\pm1}\psi_{\vep'}^*(u)^a_b\psi_\vep(u)^b_a=\delta_{\vep',\vep}.
\lb{inversionb}
\ena

Now let us consider the fusion of the vertex-face relationships \eqref{vertexface}, 
\eqref{vertexfacedual}. 
We define the $k$ fusion of the intertwining vectors by\cite{DJKMO}
 \bea
 &&\psi^{(k)}(u)^a_b=\Pi_{1\cdots k}\ \psi(u+k-1)^a_{c_1}\otimes \psi(u+k-2)^{c_1}_{c_2} \otimes \cdots 
\otimes \psi(u)^{c_{k-1}}_{b}.\lb{kpsi}
 \ena
Here the RHS is independent of the choice of $c_1,..,c_{k-1}$ provided
 $|a-c_1|=|c_1-c_2|=\cdots=|c_{k-1}-b|=1$. The local heights $a$ and $b$ now satisfy the extended admissible condition $a-b\in \{-k,-k+2,..,k\}$. 
For $k>1$, the basis $\{ v^{(k)}_\mu \}_{\mu=-k,-k+2,..,k}$ 
of $V^{(k)}$ is given by a fusion of the basis vectors $v_{\vep_i}\ (\vep_i=\pm)$ of $V$,  $k$ times. 
\bea
v^{(k)}_\mu&=&\Pi_{1\cdots k} v_{\vep_1}\otimes v_{\vep_2} \otimes \cdots \otimes v_{\vep_k}\nn\\
&=&\frac{1}{k!}\sum_{\sigma \in S_k}v_{\vep_{\sigma(1)}}\otimes v_{\vep_{\sigma(2)}}\otimes \cdots \otimes v_{\vep_{\sigma(k)}},\lb{fusionv}
\ena
where $S_k$ being the symmetric group and we set $\mu=\sum_{j=1}^k {\vep_j}$. 
Substituting \eqref{intertwinvec} to \eqref{kpsi}, we obtain
\bea
&&\psi^{(k)}(u)^a_b=\sum_{\mu\in \{-k,-k+2,..,k\}} v^{(k)}_{\mu}\psi^{(k)}_\mu(u)^a_b,
\nn\\
&&\psi^{(k)}_\mu(u)^a_b=\sum_{\vep_1,..,\vep_k=+,-\atop \mu=\sum_{j=1}^k\vep_j} 
\psi_{\vep_1}(u+k-1)^a_{c_1}\psi_{\vep_2}(u+k-2)^{c_1}_{c_2}\cdots \psi_{\vep_k}(u)^{c_{k-1}}_b.\lb{intcompo}
\ena
From \eqref{vertexface}, \eqref{kkfusionR}, 
and \eqref{kkfusionW}, it follows that
the fused intertwining vectors satisfy  
the $k\times l$ fusion vertex-face correspondence 
relations with respect to $R^{(k,l)}$ and $W^{(k,l)}$. 
% Then using \eqref{vertexface}, \eqref{hfusion}, \eqref{vhfusion}, \eqref{fhfusion}, \eqref{fvhfusion} 
% and \eqref{fusionpsi}, 
% one can verify the following.
 \bea
 &&\sum_{\mu_1',\mu_2'}
 R^{(k,l)}(u-v)_{\mu_1 \mu_2}^{\mu_1' \mu_2'}\ 
 \psi^{(k)}_{\mu_1'}(u)_{b}^{a}
 \psi^{(l)}_{\mu_2'}(v)_{c}^{b}
 =\sum_{b' \in {\mathbb{Z}}}
 \psi^{(k)}_{\mu_2}(v)_{b'}^{a}
 \psi^{(l)}_{\mu_1}(u)_{c}^{b'}
 W^{(k,l)}\left(\left.
 \begin{array}{cc}
 a&b\\
 b'&c
 \end{array}\right|u-v\right).\nn\\
 &&\lb{vertexfacekl}
 \ena

Similarly, the dual  intertwining vectors can be fused $k$ times in the following way\cite{KKW}.
 \begin{eqnarray}
\psi^{*(k)}(u)^b_a=\sum_{c_1,..,c_{k-1}}\psi^*(u+k-1)_a^{c_1}\otimes \psi^*(u+k-2)^{c_2}_{c_1} 
\otimes \cdots \otimes \psi^*(u)_{c_{k-1}}^{b}\nn\\
\lb{kpsis}
 \end{eqnarray}
 with the property 
 \be
 &&\Pi_{1 \cdots k}\ \psi^{*(k)}(u)^b_a = \psi^{*(k)}(u)^b_a\ \Pi_{1 \cdots k}. 
 \en
Written out in component form, the last relation indicates that the RHS of   
 \be
 \psi^{*(k)}_{\mu}(u)_a^b=
 \sum_{c_1,..,c_{k-1}} \psi^*_{\vep_1}(u+k-1)_a^{c_1} \psi^*_{\vep_2}(u+k-2)^{c_2}_{c_1}  \cdots  \psi^*_{\vep_k}(u)_{c_{k-1}}^{b}
 \en
 is independent of the choice of $\vep_1,\cdots,
 \vep_k$  proivided $\mu=\vep_1+\cdots+\vep_k$. 
As above, it follows immediately from \mref{vertexfacedual},
\mref{kkfusionR} and \mref{kkfusionW}, that we have
 \bea
 &&\sum_{\mu_1',\mu_2'}
 R^{(k,l)}(u-v)_{\mu_1 \mu_2}^{\mu_1' \mu_2'}\ 
 \psi^{*(k)}_{\mu_1'}(u)_{b}^{a}
 \psi^{*(l)}_{\mu_2'}(v)_{c}^{b}
 =\sum_{b' \in {\mathbb{Z}}}
 \psi^{*(k)}_{\mu_2}(v)_{b'}^{a}
 \psi^{*(l)}_{\mu_1}(u)_{c}^{b'}
 W^{(k,l)}\left(\left.
 \begin{array}{cc}
 c&b'\\
 b&a
 \end{array}\right|u-v\right).\nn\\
 &&\lb{vertexfacedualkl}
 \ena
Finally,  using \eqref{inversiona} and \eqref{inversionb}, it is easy to verify the following inversion relations.
% \begin{prop}
 \bea
 \sum_{\mu\in\{-k,-k+2,\cdots,k\}}\psi^{*(k)}_\mu(u)^a_b\psi^{(k)}_\mu(u)^b_c&=&\delta_{a,c},\lb{inversionk}
 \\
 \sum_{a\in \{b-k,b-k+2,..,b+k\}}\psi^{*(k)}_{\mu'}(u)^a_b\psi^{(k)}_\mu(u)^b_a&=&\delta_{\mu', \mu}. \lb{inversionk2}
 \ena
% \end{prop}

%For the $k=2$ case, the expressions of the intertwining vector $\psi^{(2)}(u)^a_b$ 
%and its dual $\psi^{*(2)}_\mu(u)_a^b\ (\mu=2, 0, -2)$ are given in \cite{KKW}.

\section{The Elliptic $6j$-symbols}

\subsection{The natural basis}
%We here follows Rosengren, who made a simple connection between the elliptic 
%$6j$-symbols and Spiridonov-Zhedanov's elliptic hypergeometric function. 

Let $\Theta_k$ be the space of even theta functions of order $2k$ with quasi-period $(1,\tau)$ and zero characteristics. 
\be
&&\Theta_k=\left\{f(z): {\rm entire}\ |\ f(z+1)=f(z), f(z+\tau)=e^{-2\pi i k(2z+\tau)}f(z), f(-z)=f(z)\right\}.
\en
This space forms a $k+1$ dimensional vector space. 

Let us consider the vectors of $\Theta_k$ given by
 \bea
&&e_n^k(z;\al,\beta)=[\al\pm rz]_n [\beta\pm rz]_{k-n}.\lb{natural}
%&&e_n^k(z;\al,\beta)=\prod_{j=0}^{n-1}\vtf{1}{\al+2\eta j\pm z}{\tau}
%\prod_{j=0}^{k-n-1}\vtf{1}{\beta+2\eta j\pm z}{\tau}.\lb{natural}
\ena
Through this section, we use the abbreviation 
\be
&&[u\pm z ]=[u+z][u-z],\\
&&\vartheta_\al(u\pm z|\tau)=\vartheta_\al(u+ z|\tau)\vartheta_\al(u- z|\tau),
\en
$\al=0,1,2,3$. 
The vectors $e_n^k(z;\al,\beta) \ (n=0,1,..,k)$ are linearly independent, if 
%$\eta$ is a real parameter and that 
$\al, \beta$ satisfy the following conditions\cite{Ro1}.
\be
&&\frac{\al-\beta+j}{r}\not\in \Z+\tau\Z, \ j=1-k,2-k,..., k-1, \\
&&\frac{\al+\beta+ j}{r}\not\in \Z+\tau\Z, \ j=0,1,..., k-1.
\en
Then a system of vectors $\{e_n^k(z;\al,\beta)\}_{n=0}^k$ forms a basis of $\Theta_k$, 
and is called the natural basis. 

\begin{thm}{\bf ( Elliptic binomial theorem )}\lb{Ebinomial}\cite{Ro1}

For generic parameters $\al, \beta, \gamma$, we have
\bea
&&[\al\pm rz]_k=\sum_{n=0}^k C_n^k(\al,\beta,\gamma)[\beta\pm rz]_n[\gamma\pm rz
]_{k-n} \lb{binomial}
\ena
with the coefficients
\bea
&&C_n^k(\al,\beta,\gamma)=\frac{[1]_k}{[1]_n[1]_{k-n}}\frac{[\al-\gamma,
\al+\gamma+k-n]_n[\al-\beta,\al+\beta+n]_{k-n}}{[\beta-\gamma+n-k]_n[\gamma-\beta-n]_{k-n}[\beta+\gamma]_k}.\lb{cnk}
\ena
\end{thm}
We give a proof in the Appendix.

Rosengren showed that the  change of base coefficients $R_n^m(\al,\beta,\gamma,\delta;k;q,p)$ in 
\bea
&&e^k_n(z;\al,\beta)=\sum_{m=0}^k R_n^m(\al,\beta,\gamma,\delta;k;q,p)e^k_m(z;\gamma,\delta)
\lb{changeofbase}
\ena
can be regarded as a generalisation of the elliptic $6j$-symbols. Furthermore he found an expression of the coefficients $R_n^m(\al,\beta,\gamma,\delta;k;q,p)$ 
in terms of the elliptic analogue of the very-well-poised balanced basic hypergeometric series, ${}_{12}V_{11}$. 
\begin{thm}\lb{thmRnmV}\cite{Ro1}
\bea
&&R_n^m(\al,\beta,\gamma,\delta;k;q,p)\nn\\
&&=\frac{[1]_k}{[1]_m[1]_{k-m}}\frac{[\beta-\delta,\beta+\delta-1+k]_n[\al-\gamma,\al+\gamma]_n[\beta-\gamma,\beta+\gamma]_{k-n}[\beta-\gamma]_{k-m}}{
[\gamma-\delta+m-k,\beta+\gamma]_m[\delta-\gamma-m]_{k-m}[\delta+\gamma,\beta-\gamma]_k}\nn\\
&&\times{}_{12}V_{11}(\gamma-\beta-k;-n,-m,\al-\beta+n-k,\gamma-\delta+m-k,\gamma+\delta,\al-\beta+1-k,\gamma-\beta+1). \nn\\
&&\lb{RnmV}
\ena
\end{thm}
Here ${}_{s+1}V_{s}$ is defined by\cite{Spi,SpiZhed2}
\be
&&{}_{s+1}V_{s}(u_0;u_1,\cdots,u_{s-4})=\sum_{j=0}^{\infty}\frac{[u_0+2j]}{[u_0]}\prod_{i=0}^{s-4}\frac{[u_i]_j}{[u_0+1-u_i]_j}
\en
with the balancing condition 
\be
&&\sum_{i=1}^{s-4}u_i=\frac{s-7}{2}+\frac{s-5}{2} u_0.
\en
%, which is rather convenient for dealing with
%the solvable lattice models appeared in the last section,
A proof of the Theorem is given in the Appendix. 

\subsection{Relation with the intertwining vectors}
In order to make a connection between $R_n^m(\al,\beta,\gamma,\delta;k;q,p)$ and the vertex-face intertwining vectors $\psi^{(k)}(u)^a_b$ and their duals $\psi^{*(k)}(u)^a_b$, let us consider the standard basis of $\Theta_k$ introduced by Sklyanin\cite{Sklyanin}. For $k=1$, the following two vectors form a basis of  $\Theta_1$.
\bea
&&v_+(z)=\vtf{3}{2z}{2\tau}-\vtf{2}{2z}{2\tau},\nn\\
&&v_-(z)=\vtf{3}{2z}{2\tau}+\vtf{2}{2z}{2\tau}.\lb{standard}
\ena
For $k>1$, we obtain the basis $\{ v^{(k)}_\mu(z)\}_{\mu=-k,-k+2,..,k}$ of $\Theta_k$ by fusing the basis vectors $v_\vep(z)\ (\vep=\pm)$ of $\Theta_1$. 
\bea
v^{(k)}_\mu(z)&=&\Pi_{1,2,..,k} v_{\vep_1}(z)\otimes v_{\vep_2}(z) \otimes \cdots \otimes v_{\vep_k}(z)\nn\\
&=&\frac{1}{k!}\sum_{\sigma \in S_k}v_{\vep_{\sigma(1)}}(z) v_{\vep_{\sigma(2)}}(z) \cdots v_{\vep_{\sigma(k)}}(z),\lb{fusionv}
\ena
with $\mu=\sum_{j=1}^k {\vep_j}$. 
Now let us consider the intertwining vectors in the standard basis. We set 
\bea
&&\psi(u;z)^a_b=\sum_{\vep=\pm}v_\vep(z) \psi_\vep(u)^a_b.\lb{intvecstandard}
\ena
Then from \eqref{intertwinvec} and \eqref{standard}, standard addition formulae of the theta functions yield the following formula. 
\bea
&&\psi(u;z)^a_b=\vtf{3}{z\pm\frac{(a-b)u+a}{2r}}{\tau}.\lb{psi1}
\ena

Now we consider a fusion of the intertwining vectors $\psi(u;z)^a_b$ by 
applying the same procedure as before \eqref{kpsi}. We set 
\bea
&&\psi^{(k)}(u;z)^a_b=\Pi_{1,2,..,k} \psi(u+k-1;z)^a_{c_1}\otimes\psi(u+k-2;z)^{c_1}_{c_2}\otimes \cdots \otimes \psi(u;z)^{c_{k-1}}_b.
\ena
Remember that $\psi^{(k)}(u)^a_b$ is independent of the choice of $c_1,..,c_{k-1}$ and 
$a, b$ satisfy the admissible condition $a-b\in \{-k, -k+2,..,k\}$.
Let us set $a-b=k-2n\ (n=0,1,2,..,k)$. We evaluate $\psi^{(k)}(u;z)^a_b$ in the two ways. Firstly from \eqref{intvecstandard} and \eqref{fusionv}, we obtain 
\bea
&&\psi^{(k)}(u;z)^a_b=\sum_{\mu\in \{-k,-k+2,..,k\}} v^{(k)}_\mu(z) \psi^{(k)}_\mu(u)^a_b,\lb{comp}
\ena 
where $\psi^{(k)}_\mu(u)^a_b$ is given in \eqref{intcompo}.
Secondly, from  \eqref{psi1} and a choice $c_{j+1}=c_j+1, c_0=a$ for $j=0,1,2,..,n$ and $c_{j+1}=c_j-1, c_k=b$ for $j=n, n+1, ..,k-1$, we obtain\cite{Takebe} 
\bea
\psi^{(k)}(u;z)^a_b&=&\prod_{j=0}^{n-1}\vtf{3}{\frac{-u+a-k+1+2j}{2r}\pm z}{\tau}
\prod_{j=0}^{k-n-1}\vtf{3}{\frac{u+b+2j+1}{2r}\pm z}{\tau}\nn\\
&=&(-)^ke^{-\pi i k(\frac{\tau}{2}+2z)}\prod_{j=0}^{n-1}\vtf{1}{\frac{-u+a-k+2j+1}{2r}\pm \left(z+\frac{\tau+1}{2}\right)}{\tau}\nn\\
&&\times \prod_{j=0}^{k-n-1}\vtf{1}{\frac{-u-a-k+1+2j}{2r}\pm \left(z+\frac{\tau+1}{2}\right)}{\tau}. \lb{factoredpsi}
\ena
Comparing 
%\eqref{factoredpsi} 
this with the expression \eqref{natural}, we obtain the 
identification 
\bea
\psi^{(k)}(u;z)^a_b
%&=&\sum_{\vep\in \{-k,-k+2,..,k\}} v^{(k)}_\vep(z) \psi^{(k)}_\vep(u)\nn\\
&=&(-)^ke^{-\pi i k(\frac{\tau}{2}+2z)}C^{-k} e_n^k\left(z+\frac{\tau+1}{2};\frac{-u+a-k+1}{2},\frac{-u-a-k+1}{2}\right),\nn\\
&& \lb{s2nbase}
\ena
where $C$ is a constant given in $\S$1.1.
%and $2\eta=\frac{1}{r}$.

The formulae \eqref{comp} and \eqref{s2nbase} indicate that  
 the components $\psi^{(k)}_\mu(u)^a_b$  of the vertex-face intertwining vector 
play the role of the change of base matrix elements 
from $\{v^{(k)}_\mu(z)\}$ to $\{e^k_n(z;
\frac{-u+a-k+1}{2},\frac{-u-a-k+1}{2})\}$ in $\Theta_k$. This role is similar
%, but ``twice as much'', role 
to the one of the generalized group elements (= Babelon's vertex-IRF transformations\cite{Babelon} ) 
in the theory of $q$-$6j$ symbols $\acute{\rm a}$ la Rosengren\cite{Ro}.
( See also $\S2.3$ and $\S3$ in \cite{Ro1}.)

\noindent
{\it Remark 1.} 
Using $a=b+k-2n$ in the first factors of the first line 
in \eqref{factoredpsi}, we have 
\be
\psi^{(k)}(u;z)^a_b&=&(-)^ke^{-\pi i k(\frac{\tau}{2}+2z)}\prod_{j=0}^{n-1}\vtf{1}{\frac{u-b+2j+1}{2r}\pm \left(z+\frac{\tau+1}{2}\right)}{\tau}\nn\\
&&\times \prod_{j=0}^{k-n-1}\vtf{1}{\frac{u+b+1+2j}{2r}\pm \left(z+\frac{\tau+1}{2}\right)}{\tau}. \lb{factoredpsi2}
\en
We hence obtain a second expression
\be
\psi^{(k)}(u;z)^a_b
%&=&\sum_{\vep\in \{-k,-k+2,..,k\}} v^{(k)}_\vep(z) \psi^{(k)}_\vep(u)\nn\\
&=&(-)^ke^{-\pi i k(\frac{\tau}{2}+2z)}C^{-k} e_n^k\left(z+\frac{\tau+1}{2};\frac{u-b+1}{2},\frac{u+b+1}{2}\right).
%\nn\\&& 
\lb{s2nbase2}
\en

Now let us apply the formula \eqref{s2nbase} and \eqref{comp} to \eqref{binomial}.
In the LHS, we apply \eqref{s2nbase} with $n=k$, whereas in the RHS, we
regard $\beta=\frac{-u+b-k+1}{2},\ \gamma=\frac{-u-b-k+1}{2}$ and apply 
\eqref{s2nbase} and \eqref{comp} with the replacement $a \to b,\ b \to c$. Then we get
\be
&&\psi^{(k)}_\mu(u)^a_{a+k}=\sum_{n=0}^kC^k_n(\al,\beta,\gamma)\psi^{(k)}_\mu(u)^b_c.
\en
Using the inversion relation \eqref{inversionk}, we obtain the following.
\begin{thm} 
\be
&&C^k_n(\al,\beta,\gamma)=\sum_{\mu\in \{-k,-k+2,..,k\}} \psi^{*(k)}_\mu(u)^c_b\psi_\mu^{(k)}(u)^a_{a+k},
\en
for $b-c=k-2n$ and $\al=\frac{-u+a-k+1}{2},\ \beta=\frac{-u+b-k+1}{2},\ \gamma=
\frac{-u-b-k+1}{2}$.
\end{thm}

Similarly, substituting \eqref{s2nbase} into \eqref{changeofbase} and using \eqref{comp}, we obtain  
\bea
&&\psi^{(k)}_\mu(u)^a_b=\sum_{m=0}^kR_n^m(\al,\beta,\gamma,\delta;k;q,p)
\psi^{(k)}_{\mu}(u)^c_d,
\ena
for $a-b=k-2n$, $c-d=k-2m$ and $\al=\frac{-u+a-k+1}{2}, \beta=\frac{-u-a-k+1}{2}, \gamma=\frac{-u+c-k+1}{2}, 
\delta=\frac{-u-c-k+1}{2}$. Then from the inversion relation \eqref{inversionk}, we obtain the following formula for $R_n^m$.

\begin{thm} 
\bea
&&R_n^m(\al,\beta,\gamma,\delta;k;q,p)=\sum_{\mu\in \{-k,-k+2,..,k\}} \psi^{*(k)}_\mu(u)^d_c
\psi^{(k)}_\mu(u)^a_b.\lb{Rnm}
\ena
\end{thm}

\noindent
{\it Remark 2.} This formula should be compared with the scalar product expression for $R_n^m$ derived by Rosengren ((11.2) in \cite{Ro2}), where the scaler product is defined by Sklyanin's invariant metric on $\Theta_k$. For this purpose, we need to study a scalar product formulae for the standard basis $v^{(k)}_\mu(z)$. This is an open problem.

 The expression appeared in the RHS of \eqref{Rnm} is nothing but a matrix $L^{(k)}
%\BW{a}{b}{d}{c}{u}
$ introduced by Lashkevich and Pugai\cite{LaPu} for $k=1$ and extended to higher $k$ by Kojima, Weston and the present author \cite{KKW}. Namely, 
\bea
&&L^{(k)}\BW{a}{b}{c}{d}{u}=\sum_{\mu\in \{-k,-k+2,..,k\}} \psi^{*(k)}_\mu(u)^d_c
\psi^{(k)}_\mu(u)^a_b=R_n^m(\al,\beta,\gamma,\delta;k;q,p).\nn\\
&&\lb{LRnm}
\ena
Combining \eqref{LRnm} and \eqref{RnmV}, we obtain  a full expression of 
$L^{(k)}\BW{a}{b}{c}{d}{u}$ for arbitrary $k \in \Z_{>0}$. 
\begin{cor}
For $a-b=k-2n, c-d=k-2m$, 
\be
&&L^{(k)}\BW{a}{b}{c}{d}{u}\nn\\
&&=\frac{[1]_k[-\frac{a-c}{2},-u-\frac{a+c}{2}]_m
[-\frac{a+c}{2}]_{k-m}[\frac{a-c}{2},-u+\frac{a+c}{2}-k+1]_n[-\frac{a+c}{2},-u-\frac{a-c}{2}-k+1]_{k-n}}
{[1]_m[1]_{k-m}[c+m-k,-u-\frac{a-c}{2}-k+1]_{m}[-c-m]_{k-m}[-u-k+1,-\frac{a+c}{2}]_k}
\nn\\
&&\times\sum_{j=0}^{{\rm min}(n,m) }\frac{[\frac{a+c}{2}-k+2j]}{[\frac{a+c}{2}-k]}
\frac{[\frac{a+c}{2}-k,-n,-m,a+n-k]_j}{[1,\frac{a+c}{2}+1+n-k,\frac{a+c}{2}+1+m-k,-\frac{a-c}{2}+1-n]_j}\nn\\
&&\qquad\qquad\qquad\qquad\qquad\qquad\times 
\frac{[c+m-k,-u-k+1,u+k-1,\frac{a+c}{2}+1]_j}
{[\frac{a-c}{2}+1-m,u+\frac{a+c}{2},-u+\frac{a+c}{2}+1-k,-k]_j}.
\en
\end{cor}
In \cite{KKW}, some of the $L^{(k)}$-matrix elements are calculated by fusion. 
They agree with this formula. 

\noindent
{\it Remark 3.} The $L^{(k)}$-matrix plays an important role in the 
calculation of correlation functions for the fusion eight-vertex models by using the vertex-face correspondence. 
Especially, in the limit lattice size going to infinity, a semi-infinite 
product of $L^{(k)}$ gives rise to a non-trivial operator called the tail operator 
acting on the space of states of the corner transfer matrices. Roughly speaking, 
the tail operator fills a gap between the spaces of states of the fusion eight-vertex 
and the eight-vertex SOS models. 
In addition, the  $L^{(k)}$-matrix itself is used to
write down a commutation relation between the lattice vertex operators and the 
tail operator. Based on this, the tail operator is 
identified with a certain power of the half-current of the elliptic algebra 
$U_{q,p}(\slth)$ \cite{LaPu,KKW}, 
%The algebra $U_{q,p}(\slth)$ is 
%an elliptic analogue of 
the algebra of the Drinfeld currents of 
%$U_q(\slth)$\cite{Konno}, and is 
%known to provide the Drinfeld realization of 
the face type elliptic quantum group 
${\cal B}_{q,\la}(\slth)$ \cite{JKOSII,Konno}. Such consideration 
allows us to formulate 
the lattice models by using the representation theory of $U_{q,p}(\slth)$\cite{KKW}.

Now using the formula \eqref{Rnm}, we can derive some properties of the elliptic $6j$-symbols $R_n^m(\al,\beta,\gamma,\delta;k;q,p)$. 
First of all, the inversion relation \eqref{inversionk2} imply the addition 
formula of $R_n^m$.
\begin{prop}
\bea
&&R_n^m(\al,\beta,\gamma,\delta;k;q,p)=\sum_{l=0}^k
R_n^l(\al,\beta,\rho,\sigma;k;q,p)R_l^m(\rho,\sigma,\gamma,\delta;k;q,p).
\ena
\end{prop}

Secondly, the two inversion relations \eqref{inversionk} and \eqref{inversionk2} imply the biorthogonality property of $R_n^m$ which is equivalent
to the biorthogonality condition derived first in \cite{SpiZhed}.
\begin{prop}
\bea
&&\sum_{m=0}^kR_n^m(\al,\beta,\gamma,\delta;k;q,p)R_m^l(\gamma,\delta,\al,\beta;k;q,p)
=\delta_{n,l}.
\ena
\end{prop}

Thirdly, from the fusion formulae for $\psi(u)^a_b$ and $\psi^{*}(u)^a_b$, we can verify the following fusion formula ( combinatorial formula \cite{Ro1}) for $R_n^m$.
\begin{prop}
\bea
&&R_n^m(\al,\beta,\gamma,\delta;k;q,p)\nn\\
&&=\sum_{0\leq m_j\leq 1\atop \sum_{j=1}^k m_j=m}R_{n_1}^{m_1}(\al,\al_1,\gamma,\gamma_1;1;q,p)R_{n_2}^{m_2}(\al_1,\al_2,\gamma_1,\gamma_2;1;q,p)\cdots R_{n_k}^{m_k}(\al_{k-1},\beta,\gamma_{k-1},\delta;1;q,p),\nn\\
\ena
where $n=\sum_{j=1}^k n_j\ (0\leq n_j\leq 1)$. 
\end{prop}
{\it proof.}\quad Substituting \eqref{kpsi}, 
\eqref{kpsis} to \eqref{Rnm}, the formula follows from, for example, the identification 
$\al=\frac{-u+a-k+1}{2},\ \al_{j}=\frac{-u+a_{j-1}-k+j}{2}\ (j=1,2,..,k-1),\ 
\beta=\frac{-u-a_{k-1}}{2}, \gamma=\frac{-u+c-k+1}{2},\ \gamma_j=\frac{-u+c_{j-1}-k+j}{2},(j=1,2,..,k-1),\ \delta=\frac{-u-c_{k-1}}{2} $ for $a_{j-1}-a_{j}=1-2n_j$, $c_{j-1}-c_j=1-2m_j\ (j=1,2,..,k)$ with $a_0=a, a_k=b, c_0=c, c_k=d$.\qed

Finally, using the two types of the vertex-face correspondence relation \eqref{vertexfacekl}, \eqref{vertexfacedualkl} and the inversion relations \eqref{inversionk} and \eqref{inversionk2}, we obtain the Yang-Baxter relation for $R_n^m$. 
(Figure 5)
% This is actually the face type (or dynamical) $RLL$-relation. 
\begin{thm}
\bea
&&\sum_{d}W^{(k,l)}\BW{a}{b}{d}{c}{u-v}L^{(k)}\BW{d}{c}{f}{e}{u}L^{(l)}\BW{a}{d}{g}{f}{v}\nn\\
&&=\sum_{d}L^{(k)}\BW{a}{b}{g}{d}{u}L^{(l)}\BW{b}{c}{d}{e}{v}
W^{(k,l)}\BW{g}{d}{f}{e}{u-v}.
\ena
\end{thm} 
This formula is not obvious at all from either \eqref{changeofbase} or \eqref{RnmV}.

\vspace{3mm}
\begin{center}
\setlength{\unitlength}{0.015mm}%
\begingroup\makeatletter\ifx\SetFigFont\undefined%
\gdef\SetFigFont#1#2#3#4#5{%
  \reset@font\fontsize{#1}{#2pt}%
  \fontfamily{#3}\fontseries{#4}\fontshape{#5}%
  \selectfont}%
\fi\endgroup%
\begin{picture}(8550,2581)(751,-1985)
%\put(8340,-436){\makebox(0,0)[lb]{\smash{\SetFigFont{12}{14.4}{\rmdefault}
%{\mddefault}{\updefault}{\color[rgb]{0,0,0} $\bullet$}%
%}}}
\thinlines
{\color[rgb]{0,0,0}\put(1801,239){\line( 0,-1){1200}}
}%
{\color[rgb]{0,0,0}\put(3001,239){\line( 0,-1){1200}}
}%
{\color[rgb]{0,0,0}\put(4501,239){\line( 0,-1){1200}}
}%
{\color[rgb]{0,0,0}\put(5701,239){\line( 0,-1){1200}}
}%
{\color[rgb]{0,0,0}\put(3001,239){\line(-1, 0){1200}}
}%
{\color[rgb]{0,0,0}\put(3001,-961){\line(-1, 0){1200}}
}%
{\color[rgb]{0,0,0}\put(3001,-1561){\line(-1, 0){1200}}
}%
{\color[rgb]{0,0,0}\put(5701,-961){\line(-1, 0){1200}}
}%
{\color[rgb]{0,0,0}\put(5701,239){\line(-1, 0){1200}}
}%
{\color[rgb]{0,0,0}\put(8401,-1561){\line(-1, 0){1200}}
}%
{\color[rgb]{0,0,0}\put(5701,-361){\vector(-1, 0){1200}}
}%
{\color[rgb]{0,0,0}\put(5101,239){\vector( 0,-1){1200}}
}%
{\color[rgb]{0,0,0}\put(1801,-361){\vector(-1, 0){600}}
}%
{\color[rgb]{0,0,0}\put(2401,-961){\vector( 0,-1){600}}
}%
{\color[rgb]{0,0,0}\put(8401,239){\line(-1, 0){1200}}
}%
{\color[rgb]{0,0,0}\put(7201,-361){\line( 0,-1){1200}}
}%
{\color[rgb]{0,0,0}\put(8401,-361){\line( 0,-1){1200}}
}%
{\color[rgb]{0,0,0}\put(7801,239){\vector( 0,-1){600}}
}%
{\color[rgb]{0,0,0}\put(8401,-361){\line(-1, 0){1200}}
}%
{\color[rgb]{0,0,0}\put(9001,-361){\line( 0,-1){1200}}
}%
{\color[rgb]{0,0,0}\put(9001,-961){\vector(-1, 0){600}}
}%
{\color[rgb]{0,0,0}\put(2026,239){\line(-1,-1){225}}
}%
{\color[rgb]{0,0,0}\put(7426,-361){\line(-1,-1){225}}
}%
\put(826,314){\makebox(0,0)[lb]{\smash{\SetFigFont{12}{14.4}{\rmdefault}{\mddefault}{\updefault}{\color[rgb]{0,0,0} $g$}%
}}}
\put(1651,389){\makebox(0,0)[lb]{\smash{\SetFigFont{12}{14.4}{\rmdefault}{\mddefault}{\updefault}{\color[rgb]{0,0,0} $a$}%
}}}
\put(2851,389){\makebox(0,0)[lb]{\smash{\SetFigFont{12}{14.4}{\rmdefault}{\mddefault}{\updefault}{\color[rgb]{0,0,0} $b$}%
}}}
\put(751,-1186){\makebox(0,0)[lb]{\smash{\SetFigFont{12}{14.4}{\rmdefault}{\mddefault}{\updefault}{\color[rgb]{0,0,0} $f$}%
}}}
\put(1651,-1861){\makebox(0,0)[lb]{\smash{\SetFigFont{12}{14.4}{\rmdefault}{\mddefault}{\updefault}{\color[rgb]{0,0,0} $f$}%
}}}
\put(3001,-1936){\makebox(0,0)[lb]{\smash{\SetFigFont{12}{14.4}{\rmdefault}{\mddefault}{\updefault}{\color[rgb]{0,0,0} $e$}%
}}}
\put(3151,-1186){\makebox(0,0)[lb]{\smash{\SetFigFont{12}{14.4}{\rmdefault}{\mddefault}{\updefault}{\color[rgb]{0,0,0} $c$}%
}}}
\put(4176,239){\makebox(0,0)[lb]{\smash{\SetFigFont{12}{14.4}{\rmdefault}{\mddefault}{\updefault}{\color[rgb]{0,0,0} $g$}%
}}}
\put(4501,464){\makebox(0,0)[lb]{\smash{\SetFigFont{12}{14.4}{\rmdefault}{\mddefault}{\updefault}{\color[rgb]{0,0,0} $a$}%
}}}
\put(5701,464){\makebox(0,0)[lb]{\smash{\SetFigFont{12}{14.4}{\rmdefault}{\mddefault}{\updefault}{\color[rgb]{0,0,0} $b$}%
}}}
\put(4201,-1261){\makebox(0,0)[lb]{\smash{\SetFigFont{12}{14.4}{\rmdefault}{\mddefault}{\updefault}{\color[rgb]{0,0,0} $f$}%
}}}
\put(5476,-1336){\makebox(0,0)[lb]{\smash{\SetFigFont{12}{14.4}{\rmdefault}{\mddefault}{\updefault}{\color[rgb]{0,0,0} $e$}%
}}}
\put(5926,-961){\makebox(0,0)[lb]{\smash{\SetFigFont{12}{14.4}{\rmdefault}{\mddefault}{\updefault}{\color[rgb]{0,0,0} $c$}%
}}}
\put(6901,464){\makebox(0,0)[lb]{\smash{\SetFigFont{12}{14.4}{\rmdefault}{\mddefault}{\updefault}{\color[rgb]{0,0,0} $a$}%
}}}
\put(8326,464){\makebox(0,0)[lb]{\smash{\SetFigFont{12}{14.4}{\rmdefault}{\mddefault}{\updefault}{\color[rgb]{0,0,0} $b$}%
}}}
\put(9301,-286){\makebox(0,0)[lb]{\smash{\SetFigFont{12}{14.4}{\rmdefault}{\mddefault}{\updefault}{\color[rgb]{0,0,0} $b$}%
}}}
\put(9301,-1636){\makebox(0,0)[lb]{\smash{\SetFigFont{12}{14.4}{\rmdefault}{\mddefault}{\updefault}{\color[rgb]{0,0,0} $c$}%
}}}
\put(8251,-1936){\makebox(0,0)[lb]{\smash{\SetFigFont{12}{14.4}{\rmdefault}{\mddefault}{\updefault}{\color[rgb]{0,0,0} $e$}%
}}}
\put(6826,-1936){\makebox(0,0)[lb]{\smash{\SetFigFont{12}{14.4}{\rmdefault}{\mddefault}{\updefault}{\color[rgb]{0,0,0} $f$}%
}}}
\put(3676,-661){\makebox(0,0)[lb]{\smash{\SetFigFont{12}{14.4}{\rmdefault}{\mddefault}{\updefault}{\color[rgb]{0,0,0} $=$}%
}}}
\put(6376,-661){\makebox(0,0)[lb]{\smash{\SetFigFont{12}{14.4}{\rmdefault}{\mddefault}{\updefault}{\color[rgb]{0,0,0} $=$}%
}}}
\put(6826,-211){\makebox(0,0)[lb]{\smash{\SetFigFont{12}{14.4}{\rmdefault}{\mddefault}{\updefault}{\color[rgb]{0,0,0} $g$}%
}}}
%\put(1751,-1031){\makebox(0,0)[lb]{\smash{\SetFigFont{12}{14.4}
%{\rmdefault}{\mddefault}{\updefault}{\color[rgb]{0,0,0}$\bullet$}%
%}}}
{\color[rgb]{0,0,0}\put(1201,239){\line( 0,-1){1200}}
}%
\end{picture}

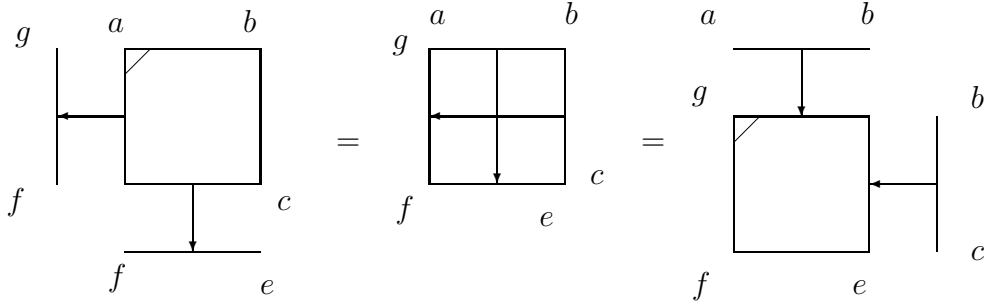
{\footnotesize Figure 5: The Yang-Baxter equation for the elliptic $6j$-symbol.}

\end{center}

\subsection{Connection with the Sklyanin algebra} 
Using the intertwining vectors and their duals, we can also construct 
a realisation of the Sklyanin algebra\cite{Bazh,Has,Qua}.  Let us define
\be
&&\cL^{(k)}(u)^b_a=\psi^{(k)}(u)^a_b\psi^{*(k)}(u)^b_a \ \in \End V^{(k)}
\en
with $a-b\in \{-k,-k+2,..,k\}$. 
Then using  \eqref{vertexfacekl}, \eqref{vertexfacedualkl} and  \eqref{inversionk},
we have
\be
&&\sum_{b}\check{R}^{(k,l)}(u-v)\cL^{(k)}(u)^b_a \otimes \cL^{(l)}(v)^c_b=\sum_b\cL^{(l)}(v)^b_a \otimes \cL^{(k)}(u)^c_b
\check{R}^{(k,l)}(u-v).
\en
Here $\check{R}^{(k,l)}(u)=P^{(k)}{R}^{(k,l)}(u)$ with $P^{(k)}$ being the permutation operator $P^{(k)} u \otimes v =v\otimes u,\  u, v\in V^{(k)}$.
Let ${\cal W}$ be a certain vector space with a basis $\{w_a \}_{a\in \Z} $. 
The  matrix $\cL^{(k)}$ can be regarded as 
a linear operator on $V^{(k)}\otimes {\cal W}$ by 
\be
&&\cL^{(k)}(u) v^{(k)}_{\mu}\otimes w_b=\sum_{\mu, a} \cL_{\mu'\mu}^{(k)}(u)^b_a \ 
v^{(k)}_{\mu'}\otimes w_a.
\en
One possible interpretation of the vector space ${\cal W}$ is given in the second paper in \cite{Has} 
as the space of theta functions spanned by the characters of the integrable representation of the affine Lie algebra $\slth$ with a fixed level. 
%We will discuss this subject elsewhere. 

%Finaly, from the inversion relation \eqref{inversionk}, it is rather obvious that 
%the intertwining vectors diagonalise the $\cL$ operators. 
%\be
%&&\cL^{(k)}(u)^b_a\psi^{(k)}(u)^a_c=\delta_{b,c} \psi^{(k)}(u)^a_c.
%\en

\vspace{3mm}
~\\
{\Large\bf Acknowledgements}~~

\noindent
The author would like to thank Hjalmar Rosengren for stimulating discussions.
He is also grateful to Masatoshi Noumi and Vyacheslav Spiridonov for their 
useful comments on the draft. He also thanks the organisers, Masatoshi Noumi 
and Kanehisa Takasaki, of the workshop ``Elliptic Integrable Systems" at Kyoto Univ., Nov. 2004, for giving him an opportunity for presenting this subject. 
There he could also make useful conversations with Jan Felipe van Diejen, Anatol Kirillov and Mikhail Olshanetsky. He is also grateful to   
Takeo Kojima and Robert Weston for their collaboration in the work \cite{KKW}. This work is supported by the Grant-in-Aid for Scientific Research (C) 15540033. 
\appendix
\setcounter{equation}{0}
\begin{appendix}
\section{Proofs of Theorem \ref{Ebinomial} and  \ref{thmRnmV}}
We here follows the idea of \cite{Ro1}. 

\subsection{Proof of Theorem \ref{Ebinomial}}
 We prove by an induction on $k$. 
Let us first split the factor $[\al+k\pm rz]$ into two parts.
\bea
[\al+k\pm rz]=A_n [\beta+n\pm rz]+B_n[\gamma+k-n \pm rz]. \lb{split}
\ena
From the theta function identity
\be
[x+y][x-y][u+v][u-v]=[u+x][u-x][v+y][v-y]-[u+y][u-y][v+x][v-x]
\en
we find 
\be
&&A_n=\frac{[\al+\gamma+k-n][\gamma-\al-n]}{[\gamma+\beta+k][\gamma-\beta+k-2n]}, \\
&&B_n= -\frac{[\beta+\al+k+n][\beta-\al+n-k]}{[\gamma+\beta+k][\gamma-\beta+k-2n]}.
\en
Then substituting \eqref{split} into 
\be
[\al\pm rz]_{k+1}
%&=&\sum_{n=0}^{k+1}C_n^{k+1}(\al,\beta,\gamma)[\beta\pm rz]_n\nn\\
&=&[\al\pm rz]_{k}[\al+k\pm rz][\gamma\pm rz]_{k-n+1}  
\en
and using \eqref{binomial}, one obtains 
\be
[\al\pm rz]_{k+1}&=&\sum_{n=0}^{k}A_{n}C^k_n(\al,\beta,\gamma)[\beta\pm rz]_{n+1}[\gamma\pm rz]_{k-n}\nn\\
&&+\sum_{n=0}^{k}B_{n}C^k_n(\al,\beta,\gamma)[\beta\pm rz]_{n}[\gamma\pm rz]_{k-n+1}.
\en
This yields the following recursion relation
\be  
&&C_n^{k+1}(\al,\beta,\gamma)=A_{n-1}C_{n-1}^k(\al,\beta,\gamma)+B_nC_n^k(\al,\beta,\gamma).
\en
Under the boundary conditions
\be
&&C_0^0=1,\ C_{-1}^k=C_{k+1}^k=0, 
\en
we find that $C_n^k(\al,\beta,\gamma)$ in \eqref{cnk} solves the recursion relation. 
\qed

%Noting that 
%\be
%&& [\al\pm rz]_k=(-)^ke^{\pi i k(\frac{\tau}{2}+2z)} C^{k}\psi^{(k)}
%(u;z-\frac{\tau+1}{2})^a_{a+k},
%\en

\subsection{Elliptic Jackson's summation formula}
\begin{thm} 
\be
&&{}_{10}V_{9}(\beta-\gamma-k;-k,\al-\gamma,-\al-\gamma+1-k,\beta+rz,\beta-rz)\nn\\
&&\qquad =\frac{[\gamma-\beta,\gamma+\beta,\al+rz,\al-rz]_k}
{[\al-\beta,\al+\beta,\gamma+rz,\gamma-rz]_k}.
\en
\end{thm}

\noindent
{\it Proof.}\quad   
Substitute $C_n^k(\al,\beta,\gamma)$ in \eqref{cnk} into 
\eqref{binomial} and compare the RHS to the definition of ${}_{10}V_9$, 
we obtain the desired formula. \hspace{\fill}$\square$

\subsection{Proof of Theorem\ref{thmRnmV}}
The equation \eqref{changeofbase} is equivalent to
\be
&&[\al\pm rz]_n[\beta\pm rz]_{k-n}=\sum_{m=0}^k 
R_n^m(\al,\beta,\gamma,\delta;k;q,p) [\gamma\pm rz]_m[\delta\pm rz]_{k-m}
\en
Substituting \eqref{binomial} to the LHS, we have
\be
&& [\al\pm rz]_n[\beta\pm rz]_{k-n}=\sum_{j=0}^nC_j^k(\al,\gamma,\beta')[\gamma\pm rz]_j[\beta'\pm rz]_{n-j}[\beta\pm rz]_{k-n}.
\en
Setting $\beta'=\beta +k-n$ and using the formula 
\be
&&[\beta\pm rz]_{l+m}=[\beta\pm rz]_{l}[\beta+l\pm rz]_{m}, 
\en
we have 
\be
&& [\al\pm rz]_n[\beta\pm rz]_{k-n}=\sum_{j=0}^nC_j^n(\al,\gamma,\beta+k-n)[\gamma\pm rz]_j[\beta\pm rz]_{k-j}.
\en
Similarly, for $[\beta\pm rz]_{k-j}$ we have 
\be
&&[\beta\pm rz]_{k-j}=\sum_{l=0}^{k-j}C_{l}^{k-j}(\beta,\delta',\delta)[\delta'\pm rz]_{l}[\delta \pm rz]_{k-j-l}.
\en
Then setting  $\delta'=\gamma+j$, we have 
\be
[\al\pm rz]_n[\beta\pm rz]_{k-n}&=&\sum_{j=0}^n \sum_{l=0}^{k-j}C_j^n(\al,\gamma,\beta+k-n)C_{l}^{k-j}(\beta,\gamma+j,\delta)[\gamma\pm rz]_{j+l}[\delta\pm rz]_{k-j-l}.
\en
Comparing this with \eqref{changeofbase}, we obtain
\be
&&R_n^m(\al,\beta,\gamma,\delta;k;q,p)=\sum_{j=0}^{{\rm min}(n,m)}C_j^n(\al,\gamma,\beta+k-n)C_{m-j}^{k-j}(\beta,\gamma+j,\delta).
\en
Substituting \eqref{cnk} into this, we obtain the formula \eqref{RnmV}.\qed

\end{appendix}

\end{document}